\documentclass[]{interact}

\usepackage{epstopdf}

\usepackage[numbers,sort&compress]{natbib}
\bibpunct[, ]{[}{]}{,}{n}{,}{,}
\makeatletter
\def\NAT@def@citea{\def\@citea{\NAT@separator}}
\makeatother

\theoremstyle{plain}
\newtheorem{theorem}{Theorem}[section]
\newtheorem{lemma}[theorem]{Lemma}

\newtheorem{proposition}[theorem]{Proposition}

\theoremstyle{definition}
\newtheorem{definition}[theorem]{Definition}

\theoremstyle{remark}
\newtheorem{remark}{Remark}

\begin{document}


\title{Calmness of the  solution mapping of Navier-Stokes  problems modeled by hemivariational inequalities}

\author{
\name{D. Inoan \textsuperscript{a}\thanks{CONTACT D. Inoan. Email: daniela.inoan@math.utcluj.ro} and J. Kolumb\'an \textsuperscript{b}}
\affil{\textsuperscript{a}Technical University of Cluj-Napoca,
             Cluj-Napoca, Romania; \textsuperscript{b}Babe\c s-Bolyai University,
              Cluj-Napoca, Romania}
}

\maketitle

\begin{abstract}
The main purpose of this paper is to find conditions for H\"older calmness of the solution mapping, viewed as a function
of the boundary data, of a hemivariational inequality governed by the Navier-Stokes operator. To this end, a more abstract model is studied first: a class of parametric equilibrium problems defined by trifunctions.
The presence of trifunctions allows the extension of the monotonicity notions and of the duality principle in the theory of equilibrium problems.
\end{abstract}


\section{Introduction}


In the papers \cite{Migorski-Ochal2005} and \cite{Migorski-Ochal2007} a class of hemivariational inequalities with the Navier-Stokes operator has been studied, where the nonslip boundary condition together with a Clarke subdifferential relation between the total pressure and the normal components of the velocity was assumed.  The main feature of such a hemivariational inequality is that it is governed by a nonmonotone and nonlinear operator, and possibly by a multivalued boundary condition defined by the Clarke derivative of a locally Lipschitz superpotential. The problem under consideration comes from fluid flow control problems and flow problems for semipermeable walls and membranes. It describes a model in which the boundary orifices in a channel are regulated to reduce the pressure of the fluid on the boundary when the normal velocity reaches a prescribed value. For the particular case when the superpotential is a lower semicontinuous convex functional the problem reduces to a variational inequality governed by a maximal monotone operator (see \cite{Chebotarev1992}, \cite{Chebotarev1997}, \cite{Konnovalova2000}).

Existence results for nonconvex locally Lipschitz superpotentials were given in the above mentioned papers \cite{Migorski-Ochal2005} (stationary case), \cite{Migorski-Ochal2007} (evolution case) and \cite{Chadli-Ansari-Yao2016}, \cite{Benaadi-Chadli-Koukkous2018}, \cite{Chadli-Kassay-Saidi2019} (periodic or antiperiodic case), for instance.

The main purpose of this paper is to study the calmness (in the sense of \cite{Clarke1976-2} and \cite{Dontchev-Rockafellar2014}) of the set of solutions of a hemivariatonal inequality governed by the Navier-Stokes operator, when the Clarke derivative is substituted by a more general control function, depending also on the state and time variables.

Calmness is an important stability property since it gives a bound on the distance between perturbed solutions and unperturbed solutions. For real functions this property is weaker than the usual local Lipschitz continuity since one of the two points considered for comparison is required to be fixed, but it is stronger than the continuity at that point. A weaker property than calmness, which we could name lower calmness, was used for the first time by F.H. Clarke in the paper \cite{Clarke1976-2}. Lower calmness is situated between calmness and lower semicontinuity. In optimization theory,  lower calmness is frequently used even in the absence of calmness (see, for instance \cite{Clarke1983}).

 In the case of set-valued mappings, calmness is defined using the excess function defined by the Romanian mathematician D. Pompeiu (see \cite{Dontchev-Rockafellar2014}). In this case calmness is a generalization of the Aubin property, which in its turn is a generalization  of the local Lipschitz property for set-valued mappings (see \cite{Dontchev-Rockafellar2014}). Calmness properties of solutions to parameterized equilibrium problems formulated with bifunctions have been studied widely. Most of this study has focused on particular models such as optimization problems (see \cite{Bonnans-Shapiro1998}, \cite{Levy-Poliquin-Rockafellar2000}, \cite{Levy2000} and \cite{Klatte1994}) and   variational or hemivariational inequalities  (see  \cite{Yen1995},  \cite {Mansour2002}).
 H\"older calmness  and H\"older continuity of solution mappings of general parametric equilibrium problems  have been studied, for instance,
in the papers \cite{Anh-Khanh2006}, \cite{Anh-Kruger-Thao2014}, \cite{Mansour-Riahi2005}, \cite{Li-Li2011}, \cite{Kassay-Bianchi-Pini2017}.  The calmness property is strongly connected to the H\"older metric subregularity of the inverse mapping  (see \cite{Dontchev-Rockafellar2014}, \cite{Zhang-Ng-Zheng-He2016} and the references therein).

In the last years many papers about hemivariational inequalites similar to that of \cite{Migorski-Ochal2005} and \cite{Migorski-Ochal2007}, with important applications, were published. As an example we mention only the quite recently appeared article \cite{Migorski-Gamorski2019}. Since the study of calmness of the solution sets for similar problems is undoubtely important, we will embed our problem in an abstract model of parametric equilibrium problem defined by  trifunction and  we will apply the abstract results obtained for the   hemivariational inequalites governed by the Navier-Stokes operator, where the parameters are functions with some properties similar to those of Clarke generalized derivatives.

Our main motivation to study equilibrium problems defined by trifunctions is the important role that the monotonicity has in existence and stability results for equilibrium problems defined by bifunctions on one hand, and the existence of bifunctions that are not monotone, on the other hand. In \cite{Inoan-Kolumban2018} we have shown that for trifunctions it is possible to define a monotonicity notion such that the monotone bifunctions that have value zero on the diagonal, generate monotone trifunctions and every bifunction is monotone as a trifunction. Therefore, for instance, the so-called mixed equilibrium problems can be formulated by monotone trifunctions. This makes possible to extend the duality principle to a large class of equilibrium problems, and to use it, for instance, by proving existence and stability of the solutions.  

This paper is structured in four sections, as follows:

 Section \ref{section:prelim} contains several notions and results needed in the sequel.

In  Section \ref{section:trifunc} we prove our main calmness result (Theorem \ref{theorem-main}), in a general setting,  for parametric equilibrium problems with trifunctions. It gives necessary and sufficient conditions for H\"older calmness based on an apriori estimation for the dual problem, which is, in fact, a gradual uniform partial calmness property (see \cite{Dontchev-Rockafellar2014}. To our best knowledge, calmness for such general problems is studied here for the first time. Existence theorems for equilibrium problems with trifunctions were given in \cite{Inoan-Kolumban2018}.

In the next section, we  apply the previous abstract theorem to mixed equilibrium problems. Finally, returning to the main purpose of the paper, we focus on the Navier-Stokes problems modeled by hemivariational type inequalities with boundary control, obtaining necessary and sufficient conditions for H\"older calmness of the solution mapping. Since  calmness is stronger than continuity, our results can be seen as sharpening of other results obtained before on the behaviour of the solutions, when the data are perturbed (see, for example, Theorem 21 in \cite{Migorski-Ochal2005} and the convergence results in \cite{Xiao-Sofonea2019}). 

\section{Preliminaries}\label{section:prelim}


For a real number $a$ we denote $a_-:=d(a,\mathbb{R}_+)=\textrm{max}\{-a,0\}$ and $a_+:=d(a,\mathbb{R}_-)=\textrm{max}\{a,0\}$.
For $a, b \in \mathbb{R}$ the following properties hold: 

(a) $(-a)_+=a_-$, $(-a)_-=a_+$;

(b) If $b\geq 0$ then $a_+\leq (a+b)_+$;

(c) $(a+b)_+ \leq a_++b_+$. The equality holds if and only if $ab\geq 0$;

(d) $(a+b)_+ \leq a_++|b|$;

(e) If $\alpha>0$, then  $a_+ \leq \alpha $ if and only if $a\leq \alpha$.

In this paper, unless otherwise mentioned, $M$ and  $X$ will be metric spaces, and - for convenience - both distances will be denoted by $d$. By $B(u,r)$ will be denoted the open ball centered at $u$, of radius $r$. The Euclidean norm on $\mathbb{R}^d$ ($d=1,2, \dots$) will be denoted by $|\cdot|$, and the scalar product of $u,v \in \mathbb{R}^d$ by $u \cdot v$.



A function $f: M\to X$ is said to satisfy the \emph{H\"older condition} of rank $k$ and exponent $\varepsilon$ if $k\geq 0$, $\varepsilon >0$, and
\begin{equation}\label{eq:holder-condition(1)}
d(f(\mu),f(\mu^\prime)) \leq k d^\varepsilon(\mu,\mu^\prime),
\end{equation}
for all $\mu,\mu^\prime \in M$. We say that $f $ is a $(k,\varepsilon)$-\emph{H\"older function near} $\bar\mu$ if there exists a neighbourhood $U(\bar \mu)$ of $\bar \mu$ such that (\ref{eq:holder-condition(1)}) is verified for all $\mu,\mu^\prime \in U(\bar \mu)$. If $\varepsilon=1$, we say that the function $f$ is $k$-Lipschitz.  A property between this and the continuity at $\bar\mu$ is the calmness at  $\bar\mu$.

A function $f: M\to X$ is said to be \emph{H\"older calm} at $\bar \mu$ if there exist $k\geq 0$, $\varepsilon >0$ and a neighborhood $U(\bar \mu)$ of $\bar \mu$ such that
\[
d(f(\mu),f(\bar\mu)) \leq k d^\varepsilon(\mu,\bar\mu),
\]
for all $\mu \in U(\bar \mu)$. If $\varepsilon=1$  we simply say that $f $  is calm instead of $(k,1)$-H\"older calm.

For a function $f:M \to \mathbb{R}$, a property strictly between calmness at $\bar\mu$ and lower semicontinuity at $\bar\mu$ is the lower calmness at $\bar\mu$. The function $f$ is said to be \emph{lower H\"older calm} at $\bar\mu$, with exponent $\varepsilon >0$, if
\begin{equation}\label{eq:lower-calm(2)}
\liminf_{\mu\to\bar\mu} \frac{f(\mu)-f(\bar\mu)}{d^\varepsilon(\mu,\bar\mu)} >-\infty.
\end{equation}
The function $f$ is said to be \emph{upper H\"older calm} at $\bar\mu$ if the function $-f$ is  lower H\"older calm at $\bar\mu$. It is clear that $f$ is H\"older calm at $\bar\mu$ if and only if it is both lower and upper H\"older calm. In optimization theory the lower calmness is frequently used even in the absence of upper calmness (see, for instance, \cite{Clarke1983}). (This could be the motivation why F.H. Clarke defined the notion of calmness by inequality (\ref{eq:lower-calm(2)}) for $M \subset \mathbb{R}$).

Calmness can be generalized for set-valued functions too.
For some sets $A,B$ in the metric space $(X,d)$ and $a \in X$, denote by
\[
d(a,B):=\inf_{b\in B}d(a,b)
\]
the distance between the  point $a$ and the set $B$, and by
\[
e(A,B):= \sup_{a\in A}d(a,B)
\]
the Pompeiu excess of $A$ with respect to $B$, where the convention
\[
e(\emptyset, B)= \left\{ \begin{array}{ll} 0, &  \textrm{when} \ B\neq \emptyset \\ +\infty, & \textrm{otherwise}  \end{array} \right.
\]
is used and $e(A, \emptyset)=+\infty$, for any set $A$, including $\emptyset$.

The Pompeiu-Hausdorff distance is defined by
\[
h(A,B)=\max (e(A,B), e(B,A)).
\]
As it is known, it does not furnish a metric on the space of all subsets, but it does on the space of nonempty closed and bounded subsets of $X$.

A set-valued mapping $S :M \to 2^X$ is said to be $(k,\varepsilon)$-\emph{H\"older continuous} if $k\geq 0$, $\varepsilon >0$, and
\[
h(S(\mu), S(\mu^\prime)) \leq k d^\varepsilon (\mu,\mu^\prime), \ \textrm{for all}\ \mu,\mu^\prime \in M.
\]
 A  mapping $S :M \to 2^X$ is said to be \emph{outer} $(k,\varepsilon)$-\emph{H\"older} at $\bar  \mu$ if $k\geq 0$, $\varepsilon >0$, and
\[
e(S(\mu), S(\bar\mu)) \leq k d^\varepsilon (\mu,\bar\mu), \ \textrm{for all}\ \mu \in M.
\]
A  mapping $S :M \to 2^X$ is said to be  $(k,\varepsilon)$-\emph{H\"older calm} at $(\bar\mu, \bar u)$ if $(\bar\mu, \bar u) \in \textrm{gph} S$, $k\geq 0$, $\varepsilon >0$, and there exist neighbourhoods $U(\bar\mu)$ of $\bar\mu$ and $V(\bar u)$ of $\bar u$  such that
\[
e(S(\mu)\cap V(\bar u), S(\bar\mu)) \leq kd^\varepsilon (\mu, \bar\mu)\ \textrm{for all} \ \mu \in U(\bar\mu).
\]
A  mapping $S :M \to 2^X$ is said to have the \emph{isolated H\"older calmness} property at $(\bar\mu, \bar u) \in \textrm{gph} S$ if there exist $k\geq 0$, $\varepsilon >0$, and the neighbourhoods $U(\bar\mu)$ of $\bar\mu$ and $V(\bar u)$ of $\bar u$  such that
\[
d(u,\bar u) \leq kd^\varepsilon (\mu, \bar\mu)\ \textrm{for all} \ \mu \in U(\bar\mu)\ \textrm{and} \ u\in S(\mu)\cap V(\bar u).
\]

If the above properties take place for $\varepsilon=1$, then instead of $(k,1)$-H\"older we say $k$-Lipschitz. For more information on these notions we propose the book \cite{Dontchev-Rockafellar2014}.

Suppose $Y$ is  a normed space and $D$ is an open subset of $Y$. If $f:D \to \mathbb{R}$ is $k$-Lipschitz near $u \in D$, then the Clarke generalized directional derivative $f^0 (u;v)$ at $u$ in the direction $v\in Y$ is defined as follows:
\[
f^0 (u;v)=\limsup_{w\to u, t \searrow 0} \frac{f(w+tv)-f(w)}{t}.
\]
Some of the basic properties of the function $f^0$ are summarised in the following
\begin{proposition}(\cite{Clarke1983}, Proposition 2.2.1 and 2.3.10)\label{prop:Clarke}
Let $f$ be $k$-Lipschitz near $u \in D$. Then

(a) The function $v \to f^0 (u;v)$ is finite, sublinear, and $k$-Lipschitz on $Y$;

(b)  The function $(u,v) \to f^0 (u;v)$ is upper semi-continuous;

(c) $f^0(u;-v)=(-f)^0(u;v)$;

(d) If $Z$ is a normed space, $T:Y\to Z$ is a surjective linear continuous operator and $g:T(D) \to \mathbb{R}$ is $k$-Lipschitz near $Tu$, then the function $g \circ T$ is $k \|T\|$-Lipschitz near $u$ and $(g\circ T)^0(u;v)=g^0(Tu,Tv)$ for $u,v\in Y$.
\end{proposition}

Usually, a bifunction $f:X \times X \to \mathbb{R}$ is said to be  \emph{monotone}  iff $f(u,v)+f(v,u) \leq 0$, for all $u,v \in X$. It is called \emph{strongly monotone} iff there exists $m >0$ such that $md^2 (u,v) \leq f(u,v)+f(v,u)$, for all $u,v \in X $. In the following definition, we extend these notions for trifunctions as follows:
\begin{definition}
The trifunction  $F:X \times X \times X  \to \mathbb{R}$ is said to be  \emph{monotone}  iff  $ F(u,v,u) \leq F(u,v,v) $ for every $u,v \in X$.

We say that $F$ is \emph{H\"older strongly  monotone}  iff there exist $m,\beta >0$ such that
\[
md^\beta (u,v) \leq F(u,v,v) - F(u,v,u), \ \textrm{for every}\ u,v \in X.
\]
\end{definition}
\begin{remark}
(a) If the trifunction $F$ has the particular form $F(u,v,w)=f(w,v)-f(w,u)$, with $f(u,u)=0$ for any $u\in X$, then $F$ is  monotone if and only if the bifunction $f$ is  monotone.

(b) If $F(u,v,w)=g(u,v)$, with $g:X \times X \to \mathbb{R}$, obviously $F$ is monotone as a trifunction. In consequence, if $G:X \times X \times X  \to \mathbb{R}$ is a monotone trifunction and $g:X \times X \to \mathbb{R}$ is arbitrary, then the trifunction $F(u,v,w)=G(u,v,w)+g(u,v) $ is also monotone. This fact simplifies, for instance, the theory of mixed equilibrium problems, and makes it more transparent (see \cite{Inoan-Kolumban2018}).

(c) If $X$ is a normed space with the dual $X^*$, then the operator $T: X \times X \to X^*$ is called \emph{semimonotone} if it is monotone with respect to the second variable, that is
\[
\langle T(u,w_1)-T(u,w_2), w_1-w_2 \rangle \geq 0, \ \textrm{for every} \ u,w_1,w_2 \in X.
\]
Variational inequalities governed by such operators were studied, for instance, in \cite{Chen1999} for single-valued functions and in \cite{Kassay-Kolumban2000} for set-valued mappings. In this last case, the function $F : X \times X \times X \to \mathbb{R}$ defined by
$ F(u,v,w)=\sup _{z \in T(u,w)}\langle z, v-u \rangle $
is monotone, while $f(u,v)=\sup _{z \in T(u,u)}\langle z, v-u \rangle$ is not monotone, so, for the variational inequality governed by $T(u,u)$, the duality principle is not applicable.
\end{remark}

\section{An abstract model}\label{section:trifunc}

In this section we consider  a general equilibrium problem  where the objective function is a trifunction that depends on a parameter $\mu$. In the papers \cite{Inoan-Kolumban2018} and \cite{Inoan-Kolumban-Minimax2018} we gave existence results for such problems, motivated by the fact that the classical theory for equilibrium problems with bifunctions can not be used for some problems that appear in applications.

 For a parameter $\mu \in M$, the problem that we study is
\[
(PE)(\mu) \quad \textrm{Find}\  \bar u \in K\ \textrm{such that} \ F(\bar u,z,\bar u; \mu) \geq 0, \
\textrm{for every} \ z \in K,
\]
where $X$ and $M$ are metric spaces,  $F:X \times X \times X \times M \to \mathbb{R}$ is a given function, and  $K$ is a nonempty subset of $X$.

Denote by $S(\mu)$  the set of solutions of the problem that depends on the parameter $\mu \in M$. Throughout the paper we suppose that it is nonempty.



\begin{theorem}\label{theorem-main} Let $\bar\mu \in M$ be a nonisolated point and $\bar u\in S(\bar \mu)$ be fixed.
Suppose that there exists a neighborhood $U(\bar \mu)$ of  $\bar\mu$, a neighborhood $\tilde V(\bar u)$ of $\bar u$ and the
 numbers $a,c, \theta \geq 0$, and $b, m,  \alpha, \beta, \xi, \theta>0$ such that

(i) $m d^\beta (\bar u,v) \leq F(\bar u,v,\bar u;\bar\mu)_-+F(\bar u,v,v;\bar\mu)_+$, for every $v \in S(\mu)\cap \tilde V(\bar u)$ and $\mu\in U(\bar \mu)$;

(ii) The estimation $\displaystyle F(\bar u,v,v;\bar\mu) \leq c d^\beta (\bar u,v)+ d^\theta (\bar u,v)[a d^\alpha (\bar u,v)+b d^{\xi}(\mu,\bar\mu)] $ holds for every $\mu\in U(\bar \mu)$ and $v \in S(\mu)\cap \tilde V(\bar u)$, with $v \neq \bar u$;

(iii)  $ \theta<\beta$ and $c<m$.

Then the mapping $S:M\to 2^X$ is H\"older calm at $(\bar\mu,\bar u)$ if and only if one of the following conditions is verified:

1) $\beta< \alpha+\theta$ and $a > 0$;

2) $\beta=\alpha+\theta$, $a>0$ and $a+c <m$;

3) $a=0$.

Moreover,  in this case we have the isolated H\"older calmness property at $(\bar \mu, \bar u)$. In the cases 2), 3) the solution $\bar u$ is unique in the neighborhood $\tilde V(\bar u)$.

 The parameters from the definition of the calmness are:

($\alpha$) $V(\bar u)=\tilde V(\bar u)$, $\displaystyle \delta=\frac{\xi}{\beta-\theta} $, $\displaystyle k=\left( \frac{b}{m-c}\right)^\frac{1}{\beta-\theta}$, for $a=0$,

($\beta$) $V(\bar u)=B(\bar u,r)\cap \tilde V(\bar u)$, $\displaystyle \delta=\frac{\xi}{\beta-\theta} $, $\displaystyle k=r \left( r^{\beta-\theta}-\frac{a}{m-c} r^\alpha \right)^\frac{1}{\theta-\beta}
\left( \frac{b}{m-c} \right)^\frac{1}{\beta-\theta}$, where $\displaystyle 0<r < \left(\frac{m-c}{a}\right)^\frac{1}{\alpha+\theta-\beta}$, for $a > 0$.
\end{theorem}
For the proof, we will need the following:


\begin{lemma}\label{prop-cuxsiy}
Let $p>0$, $q>0$, and $l\geq 0$ be given real numbers. Then the following statements are equivalent:

1) $p<q$ or $p=q$ and $l<1$, or $l=0$;

2) There exists $\varepsilon_0 >0$ such that for every $\varepsilon \in ]0, \varepsilon_0[$, there exist $\delta>0$ and $k>0$ such that, for all $x \in ]0,\varepsilon[$ and $y >0$ with $x^p-lx^q \leq y$, the inequality $x \leq ky^\delta$ holds.
\end{lemma}

\emph{Proof}: Sufficiency. 

Let $p < q$, $\varepsilon_0=l^{\frac{1}{p-q}}$ and
 let $\varphi:[0,\infty) \to \mathbb{R}$ be the function defined by $\varphi(\xi)=\xi-l \xi^\frac{q}{p}$.
It is easy to see that on the interval $]0, \varepsilon^p[$ this function has strictly positive values and is concave. From this, for $\xi \in ]0,\varepsilon^p[$, we have $\frac{\varphi(\varepsilon^p)}{\varepsilon^p} \cdot \xi <\varphi(\xi)$ that is  $\xi < \frac{\varepsilon^p}{\varphi(\varepsilon^p)} \varphi(\xi)$. Now consider $x\in ]0,\varepsilon[$ with $x^p-lx^q \leq y$  and let $\xi=x^p$. Then we get \[
x^p <\frac{\varepsilon^p}{\varphi(\varepsilon^p)} \varphi(x^p)\leq \frac{\varepsilon^p}{\varepsilon^p-l\varepsilon^q}y.
\]
and the conclusion is proved  with  \[
 k=\varepsilon (\varepsilon^p-l\varepsilon^q)^{-\frac{1}{p}}\ \textrm{and} \ \delta=\frac{1}{p}.
\]

Let $p=q$ and $l<1$. For any $x \in ]0, +\infty[$  and $y>0$, from  $x^p-lx^q \leq y$ follows $x \leq ky^\delta$, with
$k=(1-l)^{-\frac{1}{p}}$ and $\delta=\frac{1}{p}$.

Let $l=0$. For any $x \in ]0, +\infty[$  and $y>0$, from  $x^p-lx^q \leq y$ follows the conclusion with $k=1$ and $\delta=\frac{1}{p}$.

To prove the reverse implication, suppose that $p>q$. Then for $\varepsilon \in ]0, l^{\frac{1}{p-q}}[$ we have $\varphi(x^p) <0$ for all
$x \in ]0,\varepsilon[$. If 2) holds, then we would have $0<x\leq ky^\delta$, for every $y>0$, which is a contradiction. $\Box$ \\


\emph{Proof of Theorem \ref{theorem-main}}: Sufficiency. In the case $\beta< \alpha+\theta$  and $a > 0$, let $r$ be such that $0< r < (\frac{m-c}{a})^\frac{1}{\alpha+\theta-\beta}$, let the neighborhood of $\bar u$ be $V(\bar u)=B(\bar u,r)\cap \tilde V(\bar u)$, let $\mu\in U(\bar \mu)$ and $v \in S(\mu) \cap V(\bar u)$ with $v \neq \bar u$.

 Since $\bar u \in S(\bar\mu)$ and $v \in K$,
 we have
\begin{equation}\label{baru-solutie}
F(\bar u,v,\bar u;\bar\mu) \geq 0.
\end{equation}
From (i), (ii) and (\ref{baru-solutie}) follows that
\[
\begin{split}
m d^\beta(v, \bar u)& \leq F(\bar u,v,\bar u; \bar\mu)_- +
 F(\bar u, v,v; \bar\mu)_+ = F(\bar u, v,v; \bar\mu)_+ \\
 & =F(\bar u, v,v; \bar\mu) \leq c d^\beta (v, \bar u)+ d^\theta (v, \bar u)[a d^\alpha (v, \bar u)+b d^{\xi}(\mu,\bar\mu)]
 \end{split}
\]
Therefore,
\[
(m-c)d^{\beta-\theta}(v, \bar u) \leq a d^\alpha (v, \bar u)+b d^{\xi}(\mu,\bar\mu)
\]
and further on, since $m-c>0$
\begin{equation}\label{inegalitatea-de-baza}
d^{\beta-\theta}(v, \bar u) \leq \frac{a}{m-c} d^\alpha (v, \bar u)+\frac{b}{m-c} d^{\xi}(\mu,\bar\mu).
\end{equation}
We apply Lemma \ref{prop-cuxsiy}, with $x:=d(v, \bar u)$, $p:=\beta-\theta$, $q:=\alpha$, $y:=\frac{b}{m-c}d^{\xi}(\mu,\bar\mu)$, $l:=\frac{a}{m-c}$, $\varepsilon:=r$. Since $v \in B(\bar u, r)$ and $v \neq \bar u$, the conditions of the lemma are verified and we get
\begin{equation}\label{concluzia}
d(v, \bar u)\leq kd^\delta (\mu,\bar \mu),
\end{equation}
where $\displaystyle \delta=\frac{\xi}{\beta-\theta}$ and $\displaystyle k=r \left( r^{\beta-\theta}-\frac{a}{m-c} r^\alpha \right)^\frac{1}{\theta-\beta}
\left( \frac{b}{m-c} \right)^\frac{1}{\beta-\theta}$.

If $v=\bar u$, (\ref{concluzia}) is obviously verified.


In the case $\beta=\alpha+\theta$ and $a+c <m$ we have
\[
(m-c-a)d^{\beta-\theta}(v, \bar u) \leq bd^\xi (\mu,\bar \mu)
\]
which implies
\[
d(v,\bar u) \leq \left(  \frac{b}{m-c-a}\right)^\frac{1}{\beta-\theta} d^\frac{\xi}{\beta-\theta} (\mu,\bar \mu).
\]
In the case $a=0$, we can choose  $V(\bar u)= \tilde V(\bar u)$. Let $\mu\in U(\bar \mu)$ and $v \in S(\mu)\cap V(\bar u)$. In the same way as before, we get
\[
(m-c)d^{\beta-\theta}(v, \bar u) \leq b d^{\xi}(\mu,\bar\mu),
\]
so directly
\[
d(v, \bar u) \leq \left( \frac{b}{m-c}\right)^\frac{1}{\beta-\theta} d^\frac{\xi}{\beta-\theta}(\mu,\bar\mu).
\]
So, in all cases, there exist $k$ and $\delta$ such that, for every $v \in S(\mu)\cap V(\bar u)$,
\[
d(v,S(\bar \mu))=\inf_{z \in S(\bar \mu)} d(v,z) \leq d(v, \bar u)\leq kd^\delta (\mu,\bar \mu).
\]
This implies
\[
e(S(\mu)\cap V(\bar u), S(\bar\mu))= \sup_{v \in S(\mu)\cap V(\bar u)} d(v,S(\bar \mu)) \leq kd^\delta (\mu,\bar \mu).
\]
The necessity of at least one of the conditions 1), 2) or 3) follows from Lemma \ref{prop-cuxsiy} and the fact that $\bar \mu$ is nonisolated. $ \Box$
\begin{remark}
(a) If the function $F$ is H\"older strongly monotone,
then condition (i) from Theorem \ref{theorem-main} is verified. This follows directly from the fact that
$F(u,v,v) - F(u,v,u) \leq F(u,v,v)_+ + F(u,v,u)_-$. The converse is not true (see \cite{Anh-Kruger-Thao2014} for the case of bifunctions).


(b) In Section \ref{section:Navier} we will see how properties (i) and (ii) appear for hemivariational inequalities governed by the Navier-Stokes operator.

\end{remark}






\section{Parametric mixed equilibrium problems}\label{subsection-parametric-mixed}

Mixed equilibrium problems have an important role in applied mathematics. They were first studied in the paper \cite{Blum-Oettli1994}.

Consider the function $F$ having the particular form
\begin{equation}\label{kolumban2.10}
F(u,v,w;\mu)=f(w,v;\mu)-f(w,u;\mu)+g(u,v;\mu)
\end{equation}
where $f:X \times X \times M  \to \mathbb{R}$ is such that the bifunction $f(\cdot,\cdot;\mu)$ is monotone, $f(u,u;\mu)=0$, for all $u\in X $, $\mu\in M$ and $g:X \times X \times M  \to \mathbb{R}$ is an arbitrary function.  In this case we have
\begin{equation}\label{kolumban2.11}
F(u,v,u;\mu)=f(u,v;\mu)+g(u,v;\mu)\ \textrm{and} \ F(u,v,v;\mu)=-f(v,u;\mu)+g(u,v;\mu).
\end{equation}

The   problem $(PE)(\mu)$ becomes the mixed parametric equilibrium problem defined by $f$ and $g$:
\[
(PME)(\mu) \quad \textrm{Find}\  \bar u \in K\ \textrm{such that} \ f(\bar u,z;\mu)+g(\bar u,z;\mu) \geq 0, \
\textrm{for every} \ z \in K.
\]
We denote by $S(\mu)$ the set of solutions of the problem $(PME)(\mu)$ and suppose that it is nonempty.
The next result follows directly from Theorem \ref{theorem-main}.

\begin{theorem}\label{corolar-a}
Let $\bar\mu \in M$ be nonisolated and $\bar u\in S(\bar \mu)$ be fixed.
Suppose that there exist a neighborhood $U(\bar \mu)$ of  $\bar\mu$, a neighborhood $\tilde V(\bar u)$ of $\bar u$, and the
 numbers $a,b_1, b_2,c, \theta \geq 0$, $m, \alpha, \beta, \xi, \theta >0$ such that

(i) $m d^\beta (\bar u,v) \leq [f(\bar u,v;\bar \mu)+g(\bar u,v;\bar\mu)]_-+[f(v,\bar u;\bar \mu)-g(\bar u,v;\bar\mu)]_-$, for every $v \in S(\mu)\cap \tilde V(\bar u)$ and $\mu\in U(\bar \mu)$;

(ii) $\displaystyle f(v,\bar u;\mu)- f(v,\bar u;\bar \mu) \leq b_1 d^\theta (\bar u,v) d^{\xi}(\mu,\bar\mu)$, for every $\mu\in U(\bar \mu)$, $v \in S(\mu)\cap \tilde V(\bar u)$, with $v \neq \bar u$;

(iii) $\displaystyle g(\bar u,v;\bar\mu)+ g(v,\bar u;\mu) \leq c d^\beta (\bar u,v)+ d^\theta (\bar u,v)[a d^\alpha (\bar u,v)+b_2 d^{\xi}(\mu,\bar\mu)] $, for every $\mu\in U(\bar \mu)$ and $v \in S(\mu)\cap \tilde V(\bar u)$, with $v \neq \bar u$;

(iv)  $0< \beta-\theta $ and $c<m$.

Then the mapping $S:M \to 2^X$ is H\"older calm at $(\bar\mu,\bar u)$ if and only if one of the conditions 1), 2), 3) from Theorem \ref{theorem-main} is verified. Moreover, in this case we have the isolated H\"older calmness property at $(\bar \mu, \bar u)$. In the cases 2), 3) the solution $\bar u$ is unique in the neighborhood $\tilde V(\bar u)$.
\end{theorem}
\emph{Proof}:
We only have to check condition (ii) from Theorem \ref{theorem-main}.  Let $ \mu \in U(\bar \mu)$ and $v\in S(\mu)\cap \tilde V(\bar u)$. Then,
\[
 f(v,\bar u;\mu)+g(v,\bar u;\mu) \geq 0 .
 \]
We have, for $b=b_1+b_2$,
\begin{eqnarray*}
\displaystyle F(\bar u,v,v;\bar\mu) & = & -f(v,\bar u;\bar \mu)+g(\bar u,v;\bar\mu ) \\
\displaystyle  & \leq & -f(v,\bar u;\bar \mu)+g(\bar u,v;\bar\mu ) +
f(v,\bar u;\mu)+g(v,\bar u;\mu)
\\
 \displaystyle &  \leq  & c d^\beta (\bar u,v)+ d^\theta (\bar u,v)[a d^\alpha (\bar u,v)+b d^{\xi}(\mu,\bar\mu)].
\end{eqnarray*}
Therefore Theorem \ref{theorem-main} can be applied. $\Box$
\begin{remark}
Similar results on the H\"older calmness of the solution mapping have been obtained, for instance, in \cite{Anh-Kruger-Thao2014}, for  $g=0$. In this particular case, the set $K$ was considered to depend also on a parameter $\lambda$. Theorem \ref{corolar-a} can be extended in this sense, but it is not our aim in this paper.
For stability results in the case of parametric mixed problems we mention also \cite{Mansour2002}.
\end{remark}

\section{ Navier-Stokes problems modeled by hemivariational inequalities }\label{section:Navier}

In the papers \cite{Migorski-Ochal2005} and \cite{Migorski-Ochal2007},  Mig\'orski and Ochal studied a class of hemivariational problems for the Navier-Stokes operators, in the stationary and evolution case, respectively.
When $\Omega$ is  a bounded simply connected domain of $\mathbb{R}^d$, $d=2,3,\dots$, with  boundary $\Gamma$ of class $C^2$, the Navier-Stokes equations that describe the flow of a viscous incompressible constant density fluid  in the domain $\Omega$ are the following:
\begin{equation}\label{int-Navier-unu}
u^\prime-\nu \Delta u +(u\cdot \nabla)u+\nabla p=\phi,
\end{equation}
\begin{equation}\label{int-Navier-doi}
 \nabla \cdot u=0\ \textrm{on} \ Q=\Omega \times ]t_0,t_1[
\end{equation}
Here $u:\Omega \times [t_0,t_1] \to \mathbb{R}^d$ is the velocity, $\nu $ is the kinematic viscosity of the fluid, $p: \Omega \times [t_0,t_1] \to \mathbb{R}$ is the pressure, $\phi:Q \to \mathbb{R}^d$ is a vector field given by the external forces.

To obtain a variational formulation of the previous equations, it is convenient to rewrite the problem in the equivalent Leray form (see \cite{Temam1979}).

For this let us consider the set
\[
W=\{ w \in C^\infty(\Omega, \mathbb{R}^d)\ :\ \textrm{div}\, w=0\ \textrm{on} \ \Omega \}.
\]
Denote by $V$ and $H$ the closure of $W$ in the norms of $W^1_2(\Omega, \mathbb{R}^d) $ (the usual Sobolev space) and $L_2(\Omega, \mathbb{R}^d) $, respectively. We have $V \subset H \simeq H^* \subset V^*$ with the embeddings  being dense, continuous, and compact.

The space $V$ is a Hilbert space with the associated scalar product $ ((u,v))=\sum_{i=1}^d(D_iu,D_iv)$, where $D_i$ is the operator $\frac{\partial}{\partial x_i}$. Consider the spaces
\[
{\cal{V}}=L^2(t_0, t_1;V), \ {\cal{H}}=L^2(t_0, t_1;H)\ \textrm{and} \ {\cal{W}}=\{w\in {\cal{V}}\ :\ w^\prime \in {\cal{V}}^*\},
\]
where the time derivative $w^\prime$ is understood in the sense of vector valued distributions. In this case $\cal{W} \subset \cal{V} \subset \cal {H}\subset \cal{V}^*$ and the embeddings are continuous.

The pairing between $V$ and $V^*$ will be denoted by $\langle \cdot,\cdot \rangle$, and
the pairing between ${\cal{V}}$ and ${\cal{V}}^*$ will be denoted by $\langle\langle \cdot,\cdot \rangle\rangle$.
The space $\cal{W}$ is a separable reflexive Banach space with the norm
$\|w\|_{\cal{W}}=\|w\|_{\cal{V}}+\|w^\prime\|_{\cal{V}^*}$ and is continuously embedded in $C([t_0, t_1];H)$ (see \cite{Zeidler1990}).
The norm on $\mathcal{V}$ will be denoted by $\|\cdot\|$.



To write the weak formulation of the problem (\ref{int-Navier-unu})-(\ref{int-Navier-doi}), 
 consider the operators $A:V\to V^*$ and $B:V\times V \to V^*$ defined by
\begin{equation}\label{operator-A}
\langle Au, v \rangle =  \nu \int_\Omega ((u,v)) \,\textrm{d}x,
\end{equation}
\begin{equation}\label{operator-B}
\langle B(u, v),w \rangle =  \int_\Omega \sum_{i,j=1} ^d u_i(D_iv_j)w_j\, \textrm{d}x,\quad B[u]=B(u,u),
\end{equation}
and denote
\[
\langle \Phi(t), v \rangle :=\int_\Omega \phi(x,t) v(x) \textrm{d}x,
\]
for $u,v,w \in V$. It is well known (see \cite{Temam1979}, p. 162) that the operator $B$ is well defined only if $d \in \{2,3,4\}$; therefore in the following we suppose that this condition is fulfilled. For problem (\ref{int-Navier-unu})-(\ref{int-Navier-doi}) to be well posed it is necessary to assign some boundary conditions. Let us consider, for instance, in the case $d=3$, the Neumann condition $u|_\Gamma =h$, where  $h=pn-\nu \frac{\partial \gamma u}{\partial n}$, $n$ is the outward unit normal vector to $\partial \Omega$, $\frac{\partial }{\partial n}$ is the normal derivative operator, and $\gamma:V \to L^2(\Gamma)$ is the trace operator.
If we multiply the equation  (\ref{int-Navier-unu}) by a test function $v\in V$, then using the Gauss formulae, we obtain the weak formulation of the Navier-Stokes equation with the Neumann boundary condition:
\[
\langle u^\prime(t)+Au(t)+B[u(t)], v \rangle +\int_\Gamma  h \cdot \gamma v \textrm{d}\sigma(x)=\langle\Phi(t), v \rangle, \ \textrm{a.e.}\ t\in ]t_0,t_1[, v\in V.
\]
Similar to \cite{Migorski-Ochal2005} and \cite{Migorski-Ochal2007},  the Neumann boundary condition can be generalized by the subdifferential condition
\[
h(x,t)  \in \partial j(x,t,\gamma u(x,t))\ \textrm{on}\ \Gamma \times [t_0,t_1],
\]
where $\partial j $ denotes the Clarke subdifferential of the locally Lipshitz function $j : \Gamma \times [t_0,t_1] \to \mathbb{R}$. In this case the problem becomes the following evolution hemivariational inequality: For a given $\mathcal{K}\subseteq \mathcal{W}$,  find $\bar u \in \mathcal{K}$ such that
\begin{equation}\label{hemivariationala}
\langle \bar u^\prime(t)+A\bar u(t)+B[\bar u(t)]-\Phi(t), v -\bar u(t)\rangle +\int_\Gamma  j^0(x,t,\gamma \bar u(x,t); \gamma v(x)-\gamma \bar u(x,t)) \textrm{d}\sigma(x) \geq 0,
\end{equation}
for all $ v\in V$, a.e. $ t\in [t_0,t_1]$, where $j^0$ is the Clarke directional derivative of $j$. We have to note that,  although $\mathcal{K}$ is a subset of $\mathcal{W}$, in our results on hemivariational inequality (\ref{hemivariationala}), $\mathcal{K}$ will be endowed with the metric generated by the norm of $\mathcal{V}$. This is motivated by the coercivity condition verified by the operator $A$.

Using the notations (\ref{operator-A}) and (\ref{operator-B}), define the Navier-Stokes operator $N:V\to V^*$ by $Nu=Au+B[u]$,  for $u\in V$. By \cite{Temam1979}, we have the following properties:

I. $A:V\to V^*$ is linear, continuous, symmetric and $\langle Au,u \rangle \geq \nu \|u\|^2$, for all $u \in V$,

II. $B:V\times V \to V^*$ is bilinear, continuous and $\langle B(u,v), v \rangle =0$, for all $u,v \in V$,

III. The mapping $B[\cdot] :V \to V^*$ is weakly continuous.

From these follows  that the function $b$, defined by $b(u,v,z):=\langle B(u,v), z \rangle$ is trilinear and continuous. It follows that
\begin{equation}\label{eq-apare-c1}
\langle B[u]-B[v], v-u \rangle = \langle B(u-v,v), u-v \rangle \leq c_1 \cdot \|u-v\|^2\cdot \|v\|,
\end{equation}
where $c_1$ is a positive constant and $u,v \in \mathcal{V}$. We can take $c_1=\sup_{\|v\|, \|w\|=1}\langle B(w,v), w \rangle$.

 For $u,v, z\in {\cal{V}}$ we denote
\[
\langle\langle \mathcal{A} u, v \rangle\rangle =\int_{t_0}^{t_1} \langle Au(t), v(t) \rangle \textrm{d}t, \ \langle\langle \mathcal{B} (u, v), z \rangle\rangle =\int_{t_0}^{t_1} \langle B(u(t), v(t)), z(t) \rangle \textrm{d}t \ \textrm{and}
\]
\[
\langle\langle \mathcal{N} u, v \rangle\rangle =\int_{t_0}^{t_1} \langle Nu(t), v(t) \rangle \textrm{d}t.
\]

The generalized derivative $Lu=u^\prime$  defines a linear operator ${\mathcal{L}}:{\mathcal{W}} \to{\mathcal{V}}^*$ given by
\[
\langle\langle {\mathcal{L}}u,v \rangle\rangle =\int_{t_0}^{t_1} \langle u^\prime(t), v(t )\rangle \textrm{d}t,\ \textrm{for all} \ v\in {\cal{V}}.
\]

According to
\[
\langle\langle {\mathcal{L}} u, u \rangle \rangle = \int_{t_0}^{t_1}\left( \frac{1}{2} \|u(t)\|^2_H \right)^\prime \textrm{d}t =\frac{1}{2} \left( \| u(t_1)\|^2_H -\|u(t_0)\|^2_H \right)
 \]
the monotonicity on $\mathcal{K}$ of $\mathcal{L}$ follows when
for any $u_1, u_2 \in \mathcal{K}$ the inequality $\|u_2(t_0)-u_1(t_0)\|_H \leq \|u_2(t_1)-u_1(t_1)\|_H$ holds.
  This happens, for instance, in the periodic case $u(t_0)=u(t_1)$,  in the anti-periodic case $u(t_0)=-u(t_1)$, and for $\mathcal{K}=\{ u \in \mathcal{W}\ :\ u(t_0)=u_0 \}$ with a given $u_0\in H$.\\


Let $M $ be a nonempty set of functions $\mu : \Gamma \times [t_0,t_1] \times \mathbb{R}^d \times \mathbb{R}^d \to \mathbb{R}$ with the following properties:

(M1) the function $(x,t) \mapsto \mu(x,t,r;s)$ is Lebesgue measurable for all $(r,s) \in \mathbb{R}^d \times \mathbb{R}^d$;

(M2) the function $(r,s) \mapsto \mu(x,t,r;s)$ is Borel measurable  for all
$(x,t) \in \Gamma \times [t_0,t_1]$;


(M3) there exists $\theta \in [0,1]$ such that for every $\mu, \bar \mu \in M$ there exists a function $\varphi_{\mu,\bar\mu} \in L^2(\Gamma \times [t_0,t_1])$ for which
\[
|\mu(x,t,r;s)- \bar \mu(x,t,r;s)| \leq \varphi_{\mu,\bar\mu} (x,t)|s|^\theta
\]
for all $r,s \in \mathbb{R}^d $, $s$ near $0$, a.e. $(x,t) \in \Gamma \times [t_0,t_1]$.

The conditions (M1) and (M2) guarantee that the function $(x,t) \mapsto \mu (x,t, u(x,t), v(x,t))$ is Lebesque measurable for all $u \in \mathcal{H}$ (see \cite{Clarke2013}, Proposition 6.34).

For $\mu, \bar\mu \in M$ we define the distance
\begin{equation*}
d(\mu,\bar\mu) = \Big( \int_{t_0}^{t_1} \int_\Gamma
\inf_{0<\rho \leq 1} \sup_{r,s \in \mathbb{R}^d, 0<|s| \leq \rho}  |s|^{-2\theta} \cdot
 |\mu(x,t,r;s) -
  \bar \mu(x,t,r;s)|^2  \textrm{d}\sigma( x) \textrm{d}t \Big)^{\frac{1}{2}}.
\end{equation*}
From (M3) it follows that $ d(\mu,\bar\mu) < +\infty$ for all $ \mu,\bar\mu \in M$.

For $\mu\in M$ and $u,v \in {\cal{V}}$ denote
\[
\mathcal{G} _\mu (u;v)  =\int_{t_0}^{t_1} \int_\Gamma \mu (x,t,\gamma u(x,t);\gamma v(x,t) )\textrm{d}\sigma (x) \textrm{d}t.
\]

Instead of (\ref{hemivariationala}) we consider a more general problem:\\

$(NS)(\mu)$ Find $ u\in \mathcal{K}$ such that, for all $ v \in \mathcal{K}$,
\begin{equation}\label{kolumban3.1}
 \langle\langle \mathcal{L}u+ \mathcal{ N} u-\Phi, v-u \rangle\rangle + \mathcal{G} _\mu (u;v-u) \geq 0.
\end{equation}
We call such a problem  hemivariational-like inequality with boundary control variable $\mu$.
Denote by $ S(\mu)$ the set of solutions of this problem and suppose it is nonempty for all $\mu \in M$.

Using Theorem \ref{corolar-a} we are able to prove the following:


\begin{theorem} \label{teo:Navier}
Let $\mathcal{K} \subset  \mathcal{W}$ be such that the operator $\mathcal{L}$ is monotone on $\mathcal{K}$.
Let $\bar \mu \in M$ be nonisolated and $\bar u \in S(\bar \mu)$. For $\zeta >0$, define
\begin{equation}\label{eq:21}
c(\zeta):= \limsup_{v\to \bar u, v \in \mathcal{K}} \| v-\bar u\|^{-\zeta} \left( \mathcal{G}_{\bar \mu}(\bar u;v-\bar u) + \mathcal{G}_{\bar \mu}(v;\bar u-v)\right)
 \end{equation}
 and suppose that there exists $\tau >0$ such that $0<c(\tau) < +\infty$.
  Suppose further that the conditions (M1)-(M3) are verified, and there    exists $\rho >0$  such that $\|\bar u \| <\rho $ and $\rho c_1 <\nu$, where $\nu$ is the kinematic viscosity of the fluid, $c_1$ is defined in (\ref{eq-apare-c1}), and the norms $\|\varphi_{\mu,\bar\mu}\|_{L^2}$ are bounded for $\mu$ near $\bar\mu$.

  Then the mapping $\mu \mapsto S(\mu)$ is H\"older calm at $(\bar \mu, \bar u)$ if and only if one of the following conditions is verified:

1') $\tau >2$;

2') $\tau=2$ and $\rho c_1 +c(2) <\nu$;

3')   $\mathcal{G}_{\bar\mu}(\bar u,v-\bar u)+ \mathcal{G}_{\bar\mu}(v,\bar u-v) \leq 0$ for $v$ near $\bar u$.

Moreover, if one of 1'), 2') or 3') is verified, then the solution set $S$ has the isolated calmness property at $(\bar \mu,\bar u)$. In the cases 2') and 3') the solution $\bar u$ is unique.
\end{theorem}
\emph{Proof}:

Let the functions $f,g: \mathcal{K} \times \mathcal{K} \times M \to \mathbb{R}$ be defined by
\[
f(u,v):= \langle\langle \mathcal{A} u, v-u \rangle\rangle + \langle\langle {\mathcal L} u, v-u \rangle\rangle
\]
and
\[
g(u,v;\mu)= \langle\langle \mathcal{B} [u], v-u \rangle\rangle +\mathcal{G} _\mu (u;v-u) -\langle\langle \Phi, v-u \rangle\rangle.
\]
The problem is of the form (PME)($\mu$) studied in Section \ref{subsection-parametric-mixed}.

The function $f$ is strongly monotone, so  the condition (i) of Theorem \ref{corolar-a}  is verified with $\beta=2$ and $m=\nu$.
Indeed, for any $v \in \mathcal{K}$ we have
\[
\begin{split}
\nu \|\bar u-v\|^2 & \leq -f(\bar u,v)-f(v,\bar u) = -f(\bar u,v)-g(\bar u,v;\bar\mu)-f(v,\bar u)+g(\bar u,v;\bar\mu) \\
 & \leq [f(\bar u,v)+g(\bar u,v;\bar\mu)]_- +[f(v,\bar u)-g(\bar u,v;\bar\mu)]_-.
 \end{split}
\]
Condition (ii) is trivially verified.

For $\mu \in M$ and $v \in \mathcal{K}$ we have
\[
g(\bar u,v;\bar \mu) +g(v,\bar u;\mu) = \langle\langle \mathcal{B}[\bar u]-\mathcal{B}[v], v-\bar u \rangle\rangle +  \mathcal{G}_{\bar\mu}(\bar u;v-\bar u)+\mathcal{G}_\mu (v; \bar u-v).
\]
As it was mentioned before,
\[
\langle\langle \mathcal{B}[\bar u]-\mathcal{B}[v], v-\bar u \rangle\rangle = \langle\langle \mathcal{B}(\bar u-v,v), \bar u-v \rangle\rangle \leq c_1 \cdot \|\bar u-v\|^2 \cdot \|v\|.
\]
Further, by $\|\bar u \| \leq \rho$ there exists $\delta>0 $ such that, for $\|v-\bar u\| <\delta$, we have $\|v\| < \rho$.

On the other hand
\[ \begin{split}
|\mathcal{G}_{\mu}(v,\bar u-v) -\mathcal{G}_{\bar\mu}(v,\bar u-v)| & \leq \int_{t_0}^{t_1} \int_\Gamma |\mu (x,t,\gamma v(x,t); \gamma \bar u(x,t)-\gamma v(x,t)) \\
& - \bar \mu (x,t,\gamma v(x,t); \gamma \bar u(x,t)-\gamma v(x,t))| \textrm{d}\sigma(x) \textrm{d}t \\
& \leq \left( \int_{t_0}^{t_1} \int_\Gamma | \gamma \bar u(x,t)-\gamma v(x,t) |^{2\theta} \textrm{d}\sigma(x) \textrm{d}t \right)^\frac{1}{2} \cdot d(\mu,\bar\mu)\\
& \leq c_0 \|v-\bar u\|^\theta d(\mu,\bar\mu),
\end{split} \]
where $c_0=\| \gamma\|^\theta$ (see \cite{Migorski-Ochal2007} for the last inequality). Let us put
\[
a= \left \{
\begin{array}{l}
0,\ \textrm{if}\ \mathcal{G}_{\bar\mu}(\bar u,v-\bar u)+ \mathcal{G}_{\bar\mu}(v,\bar u-v) \leq 0\ \textrm{for} \ v\ \textrm{near} \ \bar u,\\
\displaystyle \frac{1}{2}(\nu-\rho c_1+c(2))\ \textrm{if}\ \tau=2\ \textrm{and} \ 0< c(2)< \nu-\rho c_1,\\
 1+c(\tau)\ \textrm{if}\ \tau=2\ \textrm{and}\ c(2) \geq \nu-\rho c_1\ \textrm{or}\ \tau \neq 2.\\
\end{array}
\right.
\]
Then, for $v\in \mathcal{K}$ near $\bar u$, we have
\[
\begin{split}
& \mathcal{G}_{\bar\mu}(\bar u;v-\bar u)+\mathcal{G}_\mu (v; \bar u-v)\\
 &  = \mathcal{G}_{\bar\mu}(\bar u;v-\bar u) +\mathcal{G}_{\bar\mu}(v, \bar u-v)-\mathcal{G}_{\bar\mu}(v, \bar u-v)+\mathcal{G}_\mu (v; \bar u-v)\\
& \leq c_0 \|v -\bar u\|^\theta \cdot d(\mu,\bar\mu) +a \|v -\bar u\|^\tau .
 \end{split}
\]
For $\|v-\bar u\| <\delta$ it follows that
\[
g(\bar u,v;\bar\mu )+ g(v,\bar u;\mu) \leq \rho c_1 \|v -\bar u\|^2 +a \|v -\bar u\|^\tau+c_0 \|v -\bar u\|^\theta \cdot d(\mu,\bar\mu).
\]
In this way, if $\|\bar u\| < \rho$ and $\|v-\bar u\| < \delta$, conditions of Theorem
\ref{corolar-a} are fulfilled with $m=\nu$, $b=c_0$, $c=\rho c_1$, $\alpha=\tau-\theta$, $\beta=2$ and $\xi=1$. $\Box$

\begin{remark}
(a) Hypothesis $\rho c_1 < \nu$ suggests that, if the viscosity coefficient $\nu$ is small, then the neighbourhood $B(0,\rho)$ of $0$, where the calmness property holds, is  small too. If $\nu$ is very small, problems may arise concerning stability and the transition towards turbulent flows (see \cite{Temam1979}). When fluctuations of flow velocity occur at very small spatial and temporal scales, one goes towards the so called turbulent models (see, for instance \cite{Berselli-Iliescu-Layton2006}, \cite{Foias-Temam2004}, \cite{Lolis-Szekelihidi2010}).

(b) Condition 3') is obviously verified if the bifunction $\mathcal{G}_{\bar \mu}$ is monotone near $\bar u$. This is why the condition $0< c(\tau) < +\infty$ from (\ref{eq:21}) is a generalization of the property which is sometimes named relaxed monotonicity (see for instance \cite{Migorski-Ochal2005}).

(c) If  the function $s \mapsto \mu(x,t,r,s)$ is positively homogeneous then
\begin{equation}\label{eq:positiv-omog}
d(\mu,\bar\mu) = \Big( \int_{t_0}^{t_1} \int_\Gamma
 \sup_{r,s \in \mathbb{R}^d, 0<|s| \leq 1}  |s|^{-2\theta } \cdot
 |\mu(x,t,r;s) -
  \bar \mu(x,t,r;s)|^2  \textrm{d}\sigma( x) \textrm{d}t \Big)^{\frac{1}{2}}
\end{equation}
(d) If the function $s \mapsto \mu(x,t,r,s)$ is linear and $\mu(x,t,r,s)=\mu_0(x,t,r)(s)$ for all $r,s \in \mathbb{R}^d$, a.e. $(x,t)\in \Gamma \times [t_0,t_1] $, then condition (M3), for $\theta=1$, can be substituted with the  following:

(M3') For every $\mu, \bar \mu \in M$, there exists a function $\varphi_{\mu,\bar\mu} \in \mathcal{V}$ such that
\[
|\mu_0(x,t,r)-\bar\mu_0(x,t,r)| \leq \varphi_{\mu,\bar\mu}(x,t),
\]
for all $r \in \mathbb{R}^d$, a.e. $(x,t)\in \Gamma \times [t_0,t_1] $. In this case we have
\[
d(\mu,\bar\mu)= \Big( \int_{t_0}^{t_1} \int_\Gamma
\sup_{r\in \mathbb{R}^d}
|\mu_0(x,t,r) -
  \bar \mu_0(x,t,r)|)^2  \textrm{d}\sigma( x) \textrm{d}t \Big)^{\frac{1}{2}},
\]
where $\bar \mu_0$ is defined similarly to $\mu_0$.

(e) If $\mu_0$ does not depend on $(x,t)$, i.e. $\mu_0(x,t,r)=\mu_1(r)$ for every $(x,t,r) \in \Gamma \times [t_0,t_1] \times \mathbb{R}^d$ then
\[
d(\mu,\bar\mu)=b_1 \sup_{r \in \mathbb{R}^d} |\mu_1(r)-\bar\mu_1(r)|,
\]
where $b_1=\textrm{mes}(\Gamma \times [t_0,t_1])$ and $\bar \mu_1$ is defined by $\bar \mu_0$ similar as $\mu_1$ by $\mu_0$.
\end{remark}

Finally, let us return to the particular case of evolution hemivariational inequalities governed by the Navier-Stokes operator. Let $J$ be a family of functions $j:\Gamma \times [t_0,t_1] \times \mathbb{R}^d  \to \mathbb{R}$ for which the function $(x,t) \mapsto j(x,t,r)$ is Lebesgue measurable for all $r \in \mathbb{R}^d$ and $r \mapsto j(x,t,r)$ is locally Lipschitz  for a.e. $(x,t)$. If we put $\mu(x,t,r;s)=j^0(x,t,r;s)$, problem $(NS)(\mu)$
reduces to (\ref{hemivariationala}). In this case, by Proposition \ref{prop:Clarke} the bifunction $\mathcal{G}_\mu : \mathcal{V} \times \mathcal{V} \to \mathbb{R}$ is well defined and, for $\theta \in [0,1]$, $r,s \in \mathbb{R}^d$, and $\mu,\bar \mu \in M$, the equality (\ref{eq:positiv-omog}) holds. As a consequence, for evolution hemivariational inequalities governed by the Navier-Stokes operator, Theorem \ref{teo:Navier}  applies.

\bibliographystyle{tfnlm}
\bibliography{bibliografia-Inoan}

\end{document}